\documentclass[11pt]{article}
\usepackage{graphicx} 
\usepackage{natbib}
\usepackage{amsmath,amsthm,amsfonts,amssymb,amscd}
\usepackage{xcolor}
\usepackage{a4wide}

\newcommand{\exponent}{{n}}

\title{
{A 5D concept for space-time optimal control problems with application to simplified Carreau flow}
}
\author{S. Beuchler, B. Endtmayer, U. Langer, A. Schafelner, T. Wick}
\date{October 30, 2025}

\begin{document}

\maketitle

\section{Status}
Optimization problems for Non-Newtonian fluids, particularly those containing nanoparticles,
arise in many applications such as in medicine, engineering, and optics. This work is related to the last topic, namely optics, which is intensively investigated in the cluster of excellence PhoenixD~\cite{Chen2018,Paetzold2017}.
Using new scientific computing techniques in optimal control we aim at the optimization of non-Newtonian fluid flows by solving the corresponding parameter dependent space-time optimality system all at once.
One model for non-Newtonian fluids is the Carreau model \cite{kennedy2005flow,ali2023carreau_roll_coating}, where we use a simplified version here.
Optimal control problems were already discussed in \cite{Beuchler2012,Endtmayer2020}, whereas partial differential equations in a space-time setting are presented in \cite{EndtmayerLangerSchafelner2024CAMWA,Endtmayer2024} by the authors. 
In this work, we consider the optimal control problem: 

Find the state $y_\rho \in Y$ and the control $u_\rho \in U$  
minimizing some cost functional
\begin{equation}
\label{Eqn:AbstractCostFunctional}
 J_\rho(y,u)
\end{equation}
subject to a simplified Carreau flow model 
\begin{equation}
\label{Eqn:StateEquation}
\begin{aligned}
     \partial_t y -\nabla_x \cdot (\left(1 + c |\nabla_x y |^2\right)^\frac{\exponent-1}{2} \nabla_x y ) &= f+u \;  \textrm{ in } Q=\Omega \times (0,T) \subset \mathbb{R}^4, \\
      y = 0 \; \textrm{ on } \partial \Omega  \times (0,T), \quad & \quad
      y(0,\cdot) =y_0=0 \;  \textrm{ in } \overline{\Omega}
\end{aligned}
\end{equation}
serves as the state equation,
where $\Omega \subset \mathbb{R}^3$ is a bounded polyhedral domain in space, $(0,T)$ denotes the time interval,
$Y$ and $U$ are appropriate spaces for the state and the control. 
We refer to \cite{Endtmayer2024,EndtmayerLangerSchafelner2024CAMWA} for the variational formulation of \eqref{Eqn:StateEquation} 
as well as existence, uniqueness and regularity results  together with the corresponding references.
The cost functional $J_\rho:Y \times U \mapsto \mathbb{R}$ depends on 
a regularization parameter $\rho$ that allows us to control the energy cost $\|u\|_U^2$ 
of the control.
The constants $n$ and $c$  are positive material parameters. 
Liquids with $\exponent>1$ are called pseudoplastic fluids,  and those  with $\exponent\leq 1$ named dilatant. 
In \cite{GuerraTiagoSequeira2014}, this was achieved for a non-simplified Carreau model on a 2D domain. 
The concept of this work consists in  the application
to a simplified Carreau model in a full space-time setting with optimization, 
and then the extension of 
the work to more sophisticated models. 
The 5D concept is structured as follows: three dimensions for space, one for time, and one for the optimization including the regularization resp. cost parameters.

\section {Current and future challenges}
Current challenges of the corresponding topic are the complexity due to nonlinearities, the optimization and the high number of unknowns arising during the simulation. Of course, a high number of unknowns increases the accuracy of the approximate solution but also comes with increased computational cost.

Using the cost functional and the state equation, we can derive 
the first-order optimality system also called  Karush-Kuhn-Tucker (KKT) system; see \cite{troeltzsch2009eng}. 
The KKT system reads as follows: Find the state $y_\rho \in Y$, the control $u_\rho \in U$ and the adjoint $p_\rho \in P$ 
such that
\begin{equation}
    \label{Eqn:KKT}
    \mathcal{L'}(y_\rho,u_\rho,p_\rho)(y,u,p)=0, \; \forall (y,u,p) \in Y \times U \times P
\end{equation}
with the Lagrangian $\mathcal{L}(y,u,p):= J_\rho(y,u)-A(y,u)(p)$,
where 
$A(\cdot,\cdot)(\cdot)$ is defined by 
\begin{equation*}
     A(y,u)(p):=(\partial_t y,p)_{L^2(Q)}+\big((1 + c |\nabla y |^2)^\frac{\exponent-1}{2} \nabla_x y, \nabla_x p\big)_{L^2(Q)}-(f+u,p)_{L^2(Q)},
 \end{equation*}
and  derived  from the variational formulation of the state equation \eqref{Eqn:StateEquation};
also see \cite{EndtmayerLangerSchafelner2024CAMWA}.

In practice, we look for an approximate solution of the KKT system  \eqref{Eqn:KKT} obtained 
by some discretization
For the discretization, here we use space-time $P_1$ finite elements, i.e. continuous, piecewise 
affine-linear simplicial finite elements
%
More precisely,
we decompose the 4D space-time cylinder $Q=\Omega \times (0,T)$ into pentatopes (4-simplices).
An example for a cost functional is the  so-called tracking type functional given by
\begin{equation*}
    J_\rho(y,u) = \frac{1}{2} \int_0^T \| y_d(t) - y(t)\|_{L^2(\Omega)}^2 dt + \frac{\rho}{2} \|u\|^2_{L^2(Q)},
\end{equation*}
where $y_d$ is the given  desired state (target), and $\rho$ denotes the regularization / cost parameter. 
In our numerical experiment, we use the desired state 
\begin{equation*}
    y_d(x,t) = \exp(-100 \|x - (0.25 * \sin(t\pi), 0.25* \cos(t\pi), t)\|_2^2),
\end{equation*}    
and the source term to $ f \equiv 0$.

\begin{figure}
    \centering
    \includegraphics[width=\linewidth]{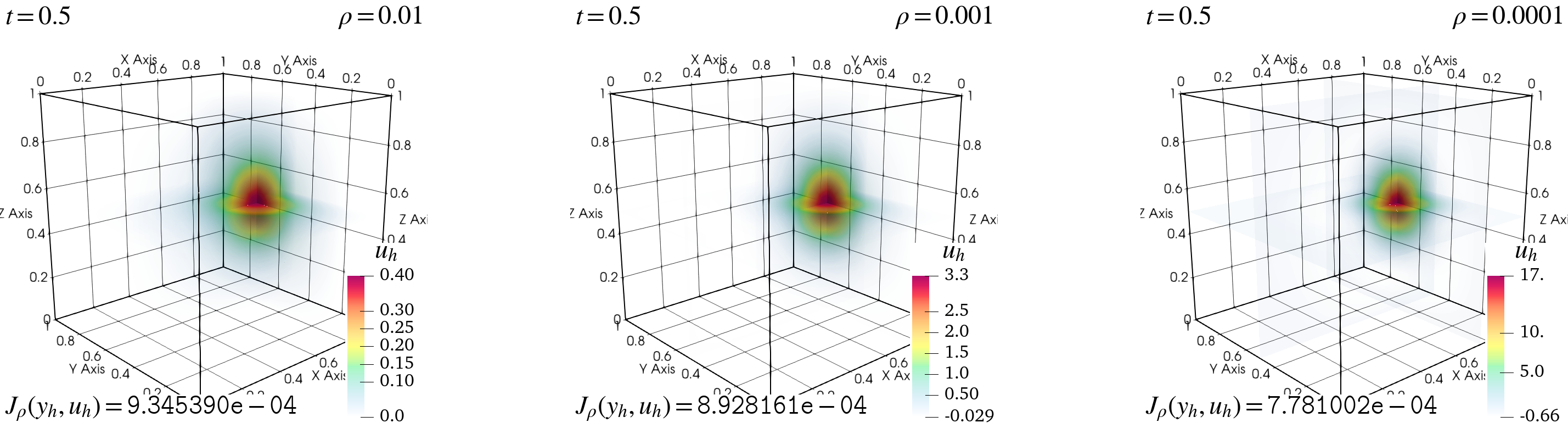}
    \caption{Optimal control $u_h$  for different values of the regularization parameter $\rho$: $\rho = 0.01$ (left); $\rho = 0.001$ (middle); $\rho = 0.0001$ (right).}
    \label{fig:placeholder}
\end{figure}
{%
In Figure~\ref{fig:placeholder}, we present a density plot of the 
control $u_h$ over the spatial domain $\Omega$ at a fixed time $t = 0.5$, for 
three points in the 5-th dimension, the regularization
parameter $\rho$, i.e. $\rho \in \{0.01,0.001,0.0001\}$. We observe a strong influence of $\rho$ 
on the magnitude of the control $u_h$. For instance, in the leftmost plot, we observe a maximal value of
$0.33$ of the control $u_h$, whereas on the rightmost plot,
the maximal value of $u_h$ has increased to $17$.
Note that the control itself has the same shape for
all $\rho$, as the target function $y_d$  is the same for all $\rho$. The simulation producing the plot in 
Figure~\ref{fig:placeholder} was done on a regular 
pentatopal grid with $\sim 480\times 10^6$ mesh elements,
resulting in approximately $31\times 10^6$ degrees of freedom (dofs) for state and costate variable, respectively. This corresponds to a global discretization parameter (mesh size) $ h \sim \frac{1}{24}$.

}
Future challenges include  
the increasing number of degrees of freedom which enhances the accuracy. However, it also comes with greater computational cost. To effectively tackle these challenges, it is crucial to develop and implement parallel algorithms capable of solving these demanding problems efficiently. Another way to improve the efficiency is adaptive mesh refinement. However, here a particular challenge is to combine it with parallel implementation, since this is not available for 4D problems until now.
Another future challenge is to use a more realistic model. However, this comes with more complex nonlinearities. This requires more research on the nonlinear solver.
Additionally, limited control problems can also be considered, where control is applied only to specific parts of the domain, such as the boundary. The practical requirements would increase the complexity of the optimization.
Moreover, additional (box) constraints imposed on the control or/and the state are of practical interest in many applications: see \cite{troeltzsch2009eng}.

\section{Advances in science and technology to meet challenges}
As computational power increases, more complex models and a higher number of degrees become more viable.

The fidelity of the model can be improved through dialogue across multiple disciplines. By including experts from various fields the model can be adapted such that it captures the real-world phenomena in a better way. However, as models become more realistic, the equations might have more nonlinear complexity, which challenges the nonlinear solver. Here, special solvers as discussed in \cite{Deuflhard2006,Cesarano2025} might be required.

 Additionally, the parallel implementation of the algorithm can significantly improve the efficiency and scalability of  the simulation. This can reduce the  computational time and therefore allow us to use the increased computational power.

\section{Concluding remarks}

In conclusion, this work presents a 5D concept to optimizing non-Newtonian fluid flows through a
simplified Carreau flow model. We solve the optimization problem by approximating the solution of the KKT System with fully space-time  finite element methods instead of the more traditional time-stepping technique combined with spatial finite element discretization. Therein, the finite element method is formulated in 3D in space, 1D in time, and 1D in the optimization loop, yielding a 5D overall framework. Finally, increased computational power, an efficient parallel
implementation and interdisciplinary collaborations will significantly advance the research in this field.

\section{Acknowledgements}

Sven Beuchler, Bernhard Endtmayer and Thomas Wick additionally acknowledge the support by the Cluster of Excellence PhoenixD (EXC 2122, Project ID390833453). 

\bibliographystyle{abbrv}
\bibliography{lit}
\end{document}